# Regulation of Disturbance Magnitude for Locational Frequency Stability Using Machine Learning

Alinane B. Kilembe
*Department of Electronic and Electrical Engineering*
*University of Strathclyde*
Glasgow, United Kingdom
alinane.brown@strath.ac.uk

Panagiotis N. Papadopoulos
*Department of Electronic and Electrical Engineering*
*University of Strathclyde*
Glasgow, United Kingdom
panagiotis.papadopoulos@strath.ac.uk

*Abstract*—Power systems must maintain the frequency within acceptable limits when subjected to a disturbance. To ensure this, the most significant credible disturbance in the system is normally used as a benchmark to allocate the Primary Frequency Response (PFR) resources. However, the overall reduction of system inertia due to increased integration of Converter Interfaced Generation (CIG) implies that systems with high penetration of CIG require more frequency control services —which are either costly or unavailable. In extreme cases of cost and scarcity, regulating the most significant disturbance magnitude can offer an efficient solution to this problem. This paper proposes a Machine Learning (ML) based technique to regulate the disturbance magnitude of the power system to comply with the frequency stability requirements i.e., Rate of Change of Frequency (RoCoF) and frequency nadir. Unlike traditional approaches which limit the disturbance magnitude by using the Centre Of Inertia (COI) because the locational frequency responses of the network are analytically hard to derive, the proposed method is able to capture such complexities using data-driven techniques. The method does not rely on the computationally intensive RMS-Time Domain Simulations (TDS), once trained offline. Consequently, by considering the locational frequency dynamics of the system, operators can identify operating conditions (OC) that fulfil frequency requirements at every monitored bus in the network, without the allocation of additional frequency control services such as inertia. The effectiveness of the proposed method is demonstrated on the modified IEEE 39 Bus network.

## I. INTRODUCTION

As the race to net-zero continues, the integration of more Converter Interfaced Generation (CIG) into the power system is forthcoming. This leads to the reduction of inertia due to the disconnection of Synchronous Generation (SG). In the Primary Frequency Response (PFR) phase, the frequency response of the network is determined by the available system inertia and the disturbance magnitude. The reduction of system inertia in the smart grid reduces the available immediate energy that resists the rapid frequency deterioration of the

This work was partially supported by the UKRI Future Leaders Fellowship MR/S034420/1 (P. N. Papadopoulos). All results can be fully reproduced using the methods and data described in this paper and the references provided. For the purpose of open access, the authors have applied for a Creative Commons Attribution (CC BY) license to any Author Accepted Manuscript version arising from this submission.

network following a disturbance [1]. In addition, it is reported that in some real networks with high integration of CIG, the frequency response dynamics have been observed as a local phenomenon rather than a global phenomenon [1]–[3].

Consequently, in a high-CIG-integrated smart grid, system operators are faced with two main challenges, i.e., i) large volumes of frequency control resources, such as inertia, are required to maintain the PFR of the system within acceptable margins, and ii) the need to consider the locational aspects of the changing system frequency dynamics instead of relying on the traditional Centre Of Inertia (COI), which poses the risk of unforeseen local protective relay activation.

To ensure the safe and reliable operation of power systems including those with high CIG penetration, [4], introduced a sequential optimisation-simulation model for optimising an operating condition (OC) with nadir considerations. In this iterative process, the method determines the minimum reserve requirement by each SG in order to adhere to frequency nadir requirements at every bus in the network. Based on the initial active power set point of the SGs, a step response is performed on every SG iteratively to accurately determine the ramp rates of the SG. The system's dynamic response is then assessed through RMS-Time Domain Simulations (TDS) within the loop. If the OC is unstable, the generator(s) dispatch is adjusted in minor increments until the minimum required ramp rate is achieved. In this work, it was effectively demonstrated that by considering the unique governor response of each SG, the approach can be economically efficient as compared to the allocation of resources based on the global minimum inertia. However, incorporating TDS within the loop can make the process computationally intensive, especially for large-scale power systems. In [5], a system optimisation model includes a set of hyperplanes describing the system frequency response requirements, i.e., nadir and Rate of Change of Frequency (RoCoF), as a function of system inertia and maximum contingency magnitude. The optimisation model is constrained using a linearised equation of the system frequency. Although this work enables the identification of OCs which, based on the system-wide response, do not violate frequency requirements without additional frequency control services, it cannot capture the locational frequency dynamics —which are increasing due to high penetration of CIG. Moreover, from an economic perspective, using the aggregated system frequency inertia without considering the location aspect on the network, may lead to higher costs as demonstrated by [3], [4].

Machine Learning (ML) models are capable of establishing complex relationships and have been used to solve different problems in power systems with a high degree of accuracy [6]– [11]. Unlike model-driven methods, ML models are data-driven approaches that do not rely on solving the Differential Algebraic Equations (DAE) of the power system to make an estimation. Consequently, ML models can make fast estimations at a very low computational budget, which arrays them as ideal candidates for real-time applications even in largescale systems, where RMS-TDS would be computationally slow [12].'

This paper proposes an ML-based approach to address the two aforementioned key challenges. The proposed method identifies OCs that do not violate frequency requirement limits at every monitored bus in the network, without the allocation of additional inertial response. This is achieved by employing two ML models; the first model captures the locational frequency dynamics of the system, while the second model captures the unique SG governor dynamic response following a disturbance. While it is analytically hard to derive the power system's locational frequency dynamics, the proposed method establishes these complex relationships using data-driven techniques. Consequently, once trained offline, the proposed approach can accurately, and with low computational burden, capture the hard-to-model frequency dynamics of the network, thereby reducing the risks associated with COI-based methods of unforeseen local relay activation.

## II. METHODOLOGY

This section presents the methodology of applying the proposed ML-based technique to regulate the maximum disturbance magnitude in a power system in order to meet the conditions for locational frequency stability. With the proposed approach, operators can identify OCs that meet the frequency stability requirements without the allocation of additional inertial response. The approach can capture the locational frequency response dynamics, elusive to COI-based methods, without the use of RMS-TDS (after offline training). This leads to fast decision-making which renders it ideal for real-time or close to real-time frequency stability monitoring even for large-scale power systems.

### A. Proposed Method Overview

The outline of the proposed three-staged methodology for ensuring conditions of locational frequency stability being met using ML-based disturbance magnitude regulation is presented in Fig. 1. Firstly, in the Initialisation Stage, the most economic OC is generated by the standard system optimisation model, such as Optimal Power Flow (OPF) —without the disturbance magnitude regulating constraints. This is then passed to the first ML model which predicts the locational frequency stability metrics of the system by estimating the nadir and RoCoF following a given disturbance. Given the fulfilment of the frequency requirements, the process is terminated, otherwise, the OC is passed to the second stage. In this stage, the second ML model predicts the unique SG governor response following the disturbance and time to nadir, i.e. governor ramp rate ($p_{rr,i}$) and ($t_{Nadir,i}$) respectively. The data from the first and second stages is then passed to the third stage which numerically estimates the maximum disturbance magnitude to ensure that there is no frequency requirement violation at any of the monitored locations in the network. This information is used to limit the disturbance magnitude in the optimisation problem by generating an updated OC. The process can be terminated after this stage, however for validation purposes, the first ML model which predicts frequency stability metrics can be re-engaged to ensure that there is no violation.

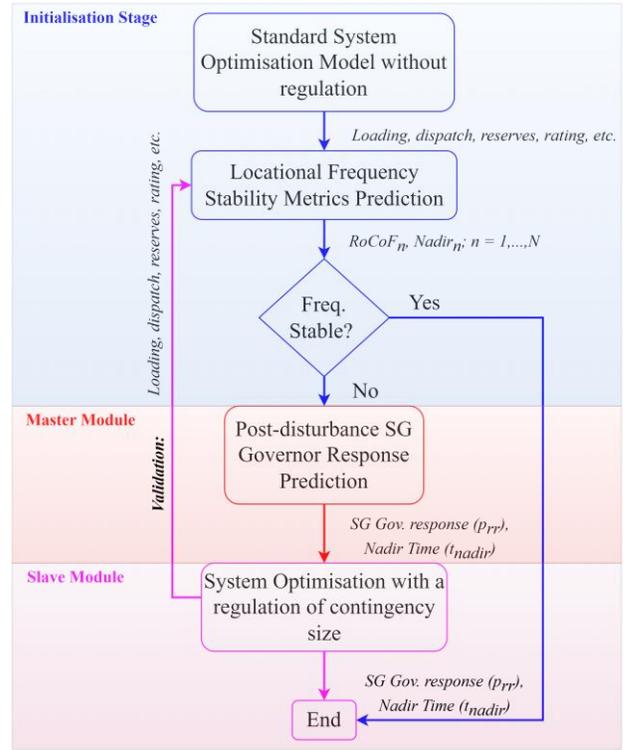

Fig. 1. A three-stage methodology for the ML-based regulation of the maximum disturbance magnitude for locational frequency stability

### B. System Frequency Response Dynamics

The power system frequency is expected to always be within the required limits following a disturbance. The rotating inertia of a power system —mainly offered by SG is the available immediate energy, $E_i$, that is injected into the power system before the activation of reserves [1]. For any given SG $i$ in the system, this is stored in the rotor and can be expressed as;

$$E_i = \frac{1}{2} J_i \omega_i^2 \qquad (1)$$

where; $J_i$ is the moment of inertia of the shaft in $kg.m^2.s$ and $\omega_i$ is the rotational speed in $rad/s$. The total amount of kinetic energy, $E_{sys}$, available at any particular moment in a power system is the sum of total rotational energy available in all the $N$ online SG i.e.,

$$E_{sys} = \sum_{i=1}^{N} \frac{1}{2} J_i \omega_i^2 \qquad (2)$$

To operate at a stable frequency, the power system is expected to maintain a balance between generation and demand. The power system dynamics of the system following a disturbance can be represented as a swing equation as follows;

$$\frac{df(t)}{dt} = \frac{P_m(t) - P_e(t)}{M_H} \qquad (3)$$

where $P_e$ represents the electric power in the system, $P_m$ represents the generator's mechanical power and $M_H$ is the inertia coefficient.

*1) RoCoF Formulation:* The frequency evolution in a network with high penetration of CIG at bus $i$, in a network with $N$ buses, following a disturbance can be stated as [2], [3];

$$\frac{df_i(t)}{dt} = \frac{df_{COI}(t)}{dt} + A_i \omega_i, \forall i = 1, ..., N \qquad (4)$$

where; $f_{COI}$ is the average frequency response from $N$ online SG, while $A_i \omega_i$ represents the oscillations or deviations, at the $i^{th}$ bus, from the COI frequency response due to distinct SG generator response in the system. In [2], it was demonstrated that it is analytically impossible to analytically derive such oscillations at every location. In this study, the ML model that predicts locational frequency stability metrics, i.e., The Locational Frequency Stability Metrics Prediction Model, effectively captures the frequency oscillations, $A\omega_i$, at every location. This factor is dependent on the system operating condition variables including; damping, distribution, locational CIG levels, generator dispatch, etc. By neglecting these oscillations while using the COI during frequency control resource allocation, unforeseen local frequency violations may occur as demonstrated in [1], [3]. Reformulating (4) to estimate individual SG responses, can be stated as;

$$\frac{df_i(t)}{dt} = \frac{\Delta P_i f_n}{2 H_i S_{ni}}, \forall i = 1, ..., N \qquad (5)$$

where $H_i$ and $S_{ni}$ are the $i^{th}$ SG inertia constant and base power in MVA respectively, while $f_n$ is the nominal system frequency. $\Delta P_i$ is the power imbalance at the generator bus $i$ following the disturbance. The $\Delta P_i$ of the $i^{th}$ SG is determined by the governor ramp rate ($p_{rr,i}$) and time to nadir, $t_{Nadir}$, derived from the Post-disturbance SG Governor Response Prediction Model. This power injected by the SG directly determines the rotor response of the generator following the disturbance. The $p_{rr,i}$ and $t_{Nadir}$ are calculated as follows;

$$p_{rr,i} = \frac{P_{Nadir,i} - P_{0,i}}{t_{Nadir,i} - t_0}, \forall i = 1, ..., N \qquad (6)$$

$$t_{Nadir_i} = \frac{\Delta P_i}{p_{rr,i}} + t_{D,i}, \forall i = 1, ..., N \qquad (7)$$

where $P_{0,i}$ is the generator output power immediately before the disturbance. $P_{Nadir,i}$ is the power output of SG $i$ at $t_{Nadir,i}$. $t_0$ is the disturbance time (which is universal), $t_{D,i}$ is the governor dead-band time. The more the SGs connected, the higher the $p_{rr}$. By monitoring and regulating $\Delta P_i$, the $i^{th}$ SG response can be limited to observe the required minimum RoCoF, $RoCoF_{limits}$, as estimated below.

$$RoCoF_{limit} = \frac{\Delta P_i - \psi_{min,i}}{\Delta P_i}.RoCoF_i, \forall i = 1, ..., N \qquad (8)$$

The $\psi_{min,i}$ is the minimum power reducing the $i^{th}$ SG's imbalance $\Delta P_i$ to observe the frequency requirements. This is the power that is redistributed to other SG in the system by the system optimisation model following a reduction of the maximum disturbance magnitude.

*2) Nadir Formulation:* The frequency nadir at the $i^{th}$ bus while incorporating the SG unique governor response can be expressed as follows;

$$\int_{t_0}^{t_{Nadir,i}} \frac{df_i(t)}{dt} = f_0 - f_{Nadir,i} \qquad (9)$$

$$f_0 - f_{Nadir,i} = \int_{t_0}^{t_{D,i}} \frac{\Delta P_i}{M_{Hi}} dt + \int_{t_{D,i}}^{t_{Nadir,i}} \frac{(p_{rr,i}.t - \Delta P_i)}{M_{Hi}} dt \qquad (10)$$

where $f_0$ is the frequency before the disturbance, while $\Delta P_i$ is the imbalance experienced by the $i^{th}$ SG. By using (6),(7),(10), the minimum nadir, $f_{Nadir,i}^{limit}$, at bus $i$ can be attained by observing the maximum imbalance, $\Delta P_i^{max} = \Delta P_i - \psi_{min,i}$, as follows;

$$f_{Nadir,i}^{limit} = f_0 - (\frac{\Delta P_i^{max}(t_0 - t_{Nadir,i}) + p_{rr,i}(t_{Nadir,i}^2 - t_{D,i})}{M_{H,i}}) \qquad (11)$$

The maximum power imbalance, $\Delta P_i^{max}$, is estimated by;

$$\Delta P_i^{max} = P_{0,i} - p_{rr,i}.(t_{Nadir,i} - \frac{f_0 - f_{Nadir,i}^{limit}}{f_0 - f_{Nadir,i}^{Old}}.t_{Nadir,i}) \qquad (12)$$

This enables the frequency nadir of the system (even at a locational level) to improve from the previous value i.e., $f_{Nadir,i}^{Old}$, to a value that is equal to or greater than the required threshold i.e., $f_{Nadir,i limit}$.

## C. Dynamic Simulations and Dataset Generation

ML models are trained and tested using a dataset of results generated from the RMS-TDS simulations. Frequency nadir, RoCoF, time to nadir (7) and SG governor response (6) following the disturbance, are recorded as regression targets of the respective ML models. System OCs are generated by varying the number of SG units committed, CIG generation and the system demand. The CIG is connected to one bus and each of the SG is an equivalent generator comprising four equal-sized units. Consequently, each SG is displaced by CIGs in four stages. The rating of $SG_{MV~Anew}$ is based on the number of remaining units, $u$, where $u \in [1,...,4]$ and is then rated to $SG_{MVA_{Aold}}$ as represented by (13). Similarly, the penetration of CIG is therefore scaled inversely as represented in (14) [13].

$$SG_{MVA,new} = u.(SG_{MVA,old}/4) \quad (13)$$

$$CIG_{MVA} = r.\frac{(5-u).SG_{MVA,old}}{4} + s.(SG_{MVA,old}) \quad (14)$$

Consequently, the penetration level is up to 40% of the overall system generation. The system loading ranges from 0.6 to 1.025 p.u. in steps of 0.25 p.u. The CIGs are all set to operate at fixed maximum active power dispatch. All the SGs maintained the default active and reactive power output limits of 0.2 to 0.85 p.u. and -0.3 to 0.7 p.u. respectively, as per rating. The initialisation of OCs was achieved by the Newton-Raphson method in DIgSILENT PowerFactory.

*D. Machine Learning Model For Capturing Frequency Dynamics*

*1) Data Pre-processing and Model Training:* In this study, two ML models have been used to 1) predict the locational frequency dynamics of the system i.e., $RoCoF_n$ and $Nadir_n$, and 2) the SG governor dynamic response i.e., ramp rate ($p_{rr,i}$) and time to nadir ($t_{Nadir,i}$). From the generated dataset, a 7030% train-test split was adopted. The number of features, $M$, used by the two ML models include the physical and steadystate characteristic variables of the network. At the time $t-1$ prior to the disturbance, the input vector x of size ($M \times 1$) is utilised by a model $\hat{p}$ to estimate the RoCoF, $\hat{R}$, and nadir, $\hat{N}$, of the $N$ buses of the network as shown below;

$$\hat{R}_n, \hat{N}_n = \hat{p}(l_{t-1}, g_{t-1}, d_{t-1}, o_{t-1}), \forall n = 1,...,N \quad (15)$$

where; $l$ is system loading, $g$ is CIG output, $d$ is SG dispatch, and $o$ are generator ratings. Assessment of model performance was conducted using $Cross-Validation$ (5-fold) which randomly splits the training set into $k-fold$, whereby each $k-1$ set is used for training, with the remaining set used for testing [14]. The $sklearn-GridSearchCV$ function [14], is used to optimise the performance of each ML algorithm. Standardisation of the dataset is achieved by using the $Standard-Scaler$ function which scales every variable to unit variance. The stored mean and standard deviation, through Inverse-Transform, are used to re-scale the data, for testing and evaluation. The scaling process is represented as;

$$z_i = s_i.\frac{x_i - \mu_i}{\sigma_i}, \forall i \in \{1,...,M\} \quad (16)$$

where; $z_i$ is the standard score, $x_i$ is the value, $\mu_i$ is the mean, $\sigma_i$ is the standard deviation and $s_i$ is the scaling factor of the $i^{th}$ feature or regression target. $M$ represents the total number of dataset features and targets.

*2) Multilayer Perceptron (MLP):* The proposed method utilises two MLP neural networks to capture locational

frequency dynamics and unique SG governor response. This data is used to determine the maximum disturbance magnitude of the system to satisfy the conditions for locational frequency stability. MLP is a feed-forward Artificial Neural Network (ANN) that contains a minimum of three layers; the input layer, the hidden layer with $\sigma_v$ neurons each, and the output layer. For an MLP model with hidden layers, $v$, feature vector x, weights matrix $W_v$ of size ($\sigma_v \times \sigma_{v+1}$) and bias vector $b_v$ of size ($\sigma_v \times 1$), it can be represented by (17-19).

$$z_1 = \mathbf{W_1}^T\mathbf{x} + \mathbf{b_1} \quad (17)$$

$$\hat{z}_{v+1} = \mathbf{W_{v+1}}^T\mathbf{z}_v + \mathbf{b}_{v+1}, \forall v = 1,..,V-1 \quad (18)$$

The Rectified Linear Unit (ReLU) activation function, $\Theta$, which easily overcomes numerical problems associated with the sigmoid is chosen. ReLU is expressed as $\Theta(\hat{z}_v) = max(0, \hat{z}_v), \forall v = 1,...,V-1$. The predicted frequency stability metrics ($\hat{N}, \hat{R}$) vector, y ($1 \times 2N$), is therefore given as;

$$\mathbf{y} = \mathbf{W}^T_{V+1}\mathbf{z}_V + \mathbf{b}_{V+1}, \mathbf{y} = [\hat{R}_n, \hat{N}_n], \forall n = 1,...,N \quad (19)$$

The first step in MLP is to propagate the features up to the output layer, a step known as forward propagation. Thereafter, based on the output, an error, $E$, is calculated whereby the weights in the hidden layer(s) are adjusted to minimise the same (back-propagation) [14]. This is achieved through the calculation of the error derivative of each weight at a specified learning rate, $\eta$, in (20).

$$w_i \leftarrow w_i - \eta \nabla E \quad (20)$$

Finally, these steps are repeated several times over epochs to establish the best model parameters with two outputs, i.e., *RoCoF* and *nadir*. In the study, the fully connected MLP architecture used had three hidden layers with 100 *neurons* each, 0.001 *alpha*, 0.01 *learningrate* and a maximum iteration of 2000 which was determined during $sklearn-GridSearchCV$ hyperparameter tuning.

*E. Model Accuracy Evaluation*

The performance of the two ML models is evaluated using the Root Mean Squared Error (RMSE) metric. For $J$ OCs in the testing dataset, the RMSE between the actual variable ($y$) and the predicted variable ($\hat{y}$) is given by (21).

$$RMSE_i = \sqrt{\frac{\sum_{n=1}^{J}(y_{n,i} - \hat{y}_{n,i})^2}{J}} \quad (21)$$

Where; $J$ is the total number of OCs in the dataset and $i$ is the regression target. Errors in critical cases, i.e., those close to the stability boundary, impact key decisions, such as ancillary service procurement. An overestimate of RoCoF and/or nadir may result in the operator over-procuring costly ancillary services. Conversely, an underestimate may result in the operator procuring insufficient ancillary services, leaving the system potentially vulnerable.

## III. Results

### A. Case Study

AC OPF was implemented in MATPOWER (using MIPS algorithm) and RMS-TDS were conducted in DIgSILENT on a modified IEEE 39-bus network as a test bed [13]. The penetration of CIG is achieved through the connection at bus 16 in the highlighted Area 2 of Fig.2. The Western Electricity Coordinating Council (WECC) Type 4 Wind Turbine Generator control model [7] is used to connect the CIG to the grid through a full converter interface.

Each of the SG [4, 5, 6, 7] is displaced by the CIG in four levels as shown by (13) and (14). The instantaneous penetration of CIG ranges from 100 MW to 1000 MW, representing close to 40% of the system generation. The generator outage event is SG 5 which makes a significant generation contribution in the region of up to 25%. The thresholds for RoCoF protection and protective Under Frequency Load Shedding (UFLS) are considered to be -0.5 Hz/s and 59.6 Hz respectively [15]. All simulations are carried out on an 11th Gen Intel (R) Core (TM) i7-11700 @ 2.50 GHz with 16 GB installed RAM, which took close to 6.5 hours to execute 2,200 simulations and generate the database, extract key features and targets for the training and testing of the two ML models. The Python programming language was used for data preprocessing, extracting useful features of interest, training the ML algorithms, and data analysis.

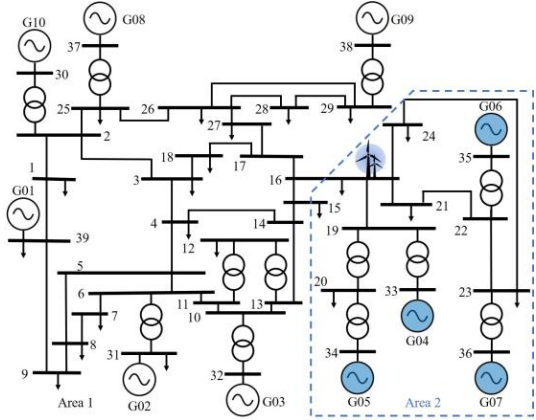

Fig. 2. The modified IEEE 39-Bus Network highlighting the region with the CIG location at Bus 16

### B. Accuracy of the ML Models in Predicting Frequency Dynamics

The RMSE is used to evaluate the accuracy of the two ML models in predicting the locational frequency metrics of the system i.e., *RoCoF* and *Nadir* and the SG governor dynamic response i.e., governor *RampRate* and *TimeToNadir*. Accurate prediction of the locational frequency stability metrics ensures the correct estimation of the network disturbance/imbalance propagation, $\Delta P_i$. This information is essential to determine the permissible maximum disturbance that ensures no locational frequency limit violations. In Table I it

| Bus | RoCoF(Hz/s) | Nadir(Hz) | Nadir time(s) | Ramping(MW/s) |
|-----|-------------|-----------|---------------|---------------|
| 30 | 0.0009 | 0.0078 | 0.2636 | 0.0480 |
| 31 | 0.0008 | 0.0084 | 0.2714 | 0.3196 |
| 32 | 0.0009 | 0.0078 | 0.2788 | 0.1547 |
| 33 | 0.0014 | 0.0083 | 0.2817 | 0.0305 |
| 34 | 0.0012 | 0.0078 | 0.2372 | 0.4817 |
| 35 | 0.0013 | 0.0078 | 0.2798 | 1.2844 |
| 36 | 0.0014 | 0.0077 | 0.2751 | 0.1364 |
| 37 | 0.0009 | 0.0082 | 0.2738 | 0.0761 |
| 38 | 0.0011 | 0.0077 | 0.3037 | 0.0389 |
| 39 | 0.0008 | 0.0078 | 0.2602 | 0.0000 |

is observed that the Locational Frequency Stability Metrics Prediction Model has the maximum RoCoF RMSE of 0.0014 Hz/s at Bus 33 and 36. In addition, the maximum RMSE for the frequency nadir metric is 0.0084 Hz observed at Bus 31, thus, demonstrating a high degree of accuracy by this model.

On the other hand, the Post-disturbance SG Governor Response Prediction Model, which predicts *TimeToNadir* and SG *RampRate*, portrays a slightly lower accuracy. Over a simulation window of 60 seconds, the model has the maximum *TimeToNadir* RMSE of 0.3037 seconds at Bus 38, representing 0.5062% error. Similarly, the maximum *RampRate* RMSE is 1.2844 MW/s at Bus 35 whose SG is rated 800 *MV A*, representing 0.1606% error. The accuracy of this model directly affects the accuracy of the estimated maximum disturbance magnitude. This accuracy can be improved in several ways including; utilising separate models in predicting *RampRate* and *TimeToNadir*, and using separate models at each location as it was done in [13] for transient stability.

Concerning the computational performance for online application in frequency stability metrics prediction, a sample of 500 OCs is used. It is found that 300 seconds is required by the RMS-TDS while 0.0090 seconds is required by the ANN. As a data-driven technique, the ANN does not need to solve the DAEs of the network to make an estimation. Consequently, it is fast at making estimations even at low computation budget, hence, a good candidate for online applications.

### C. Regulation of the Maximum Disturbance Magnitude for Frequency Stability

RMS-TDS are conducted to evaluate the performance of the proposed method outlined in Fig. 1. OCs with frequency violation are sampled from the test dataset and used to evaluate the dynamic response of the system following the adjusted maximum disturbance. By implementing the proposed method, i.e., the Regulated Model, results from the OPF are exported to DIgSILENT for validation. It can be observed in Table II that the model without disturbance magnitude regulation, i.e., the Unregulated Model, has a minimum nadir of 59.37 Hz and maximum nadir of 59.49 Hz —representing 100% of frequency requirement violations. These violations have the potential to activate the UFLS relays and lead to large-scale blackouts.

Contrariwise, after implementing the Regulated Model, better results are achieved with an accuracy improving

TABLE II

ACCURACY OF THE PROPOSED REGULATED SYSTEM OPTIMISATION MODEL

| Assessment Criteria | Value |
|---|---|
| Nadir Boundary (Hz) | 59.60 |
| Regulated Model Minimum Nadir (Hz) | 59.59 |
| Unregulated Model Minimum Nadir (Hz) | 59.37 |
| Regulated Model Maximum Nadir (Hz) | 59.77 |
| Unregulated Model Maximum Nadir (Hz) | 59.49 |
| Maximum RoCoF estimation Error (Hz/s) | 0.05 |
| Minimum RoCoF estimation Error (Hz/s) | $1.12 \times 10^{-4}$ |
| Regulated Model Accuracy (%) | 98.77 |

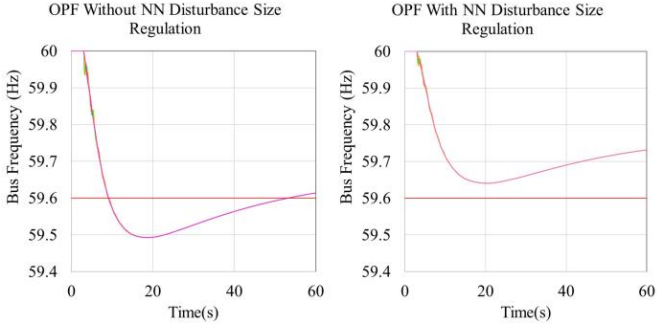

Fig. 3. Comparison of locational frequency response (No disturbance magnitude regulation (left), With disturbance regulation (right)) following an outage of SG 05 with CIG generation of 617 MW at 0.8 system loading

from 0% up to 98.77%. Out of the 1.23% that caused a violation, the minimum nadir is 59.59 Hz, which is 0.01 Hz within the stability limits. The errors portrayed by the Regulated Model are attributed to the combined predictive error of the two ML models. This error can be improved by slightly overestimating the threshold as done by [3] and [7] for locational frequency and transient stability problems respectively. An example of the performance of the two models is given in Fig. 2. It is observed that the Regulated Model fulfils conditions for locational frequency stability, unlike the Unregulated Model. This demonstrates that the former can effectively identify a feasible OC without the allocation of additional frequency control resources such as inertia.

## IV. CONCLUSION

In a system with high penetration of CIG, the requirements of Primary Frequency Response (PFR) —which are costly and/or scarce, can be very significant. In this study, we propose an ML technique-based method to identify stable operating conditions (OCs), concerning frequency stability, without the allocation of additional inertia response. The method can fulfil the conditions for locational frequency stability by capturing the locational frequency dynamics. This approach captures and considers the detailed frequency dynamics of the network in the estimation of the maximum disturbance magnitude, which is overlooked by COI-based methods, thereby bearing the risk of unforeseen local relay activation. Two separate Artificial Neural Networks (ANN) are used to predict the locational frequency stability metrics, i.e., $RoCoF_n$ and $Nadir_n$, as well as the governor response of online SG, i.e. $RampRate_i$ and $TimeToNadir_i$. By monitoring and regulating the maximum power imbalance at each bus bar connected to SG $i$, the maximum magnitude of the disturbance necessary to observe the frequency requirements at the locational level is determined. Thereafter, the system optimisation model (with disturbance regulation, i.e., Regulated Model) is used to re-dispatch the generators. The proposed method does not rely on RMSTDS after offline training of the ML models. This makes it faster and ideal for real-time and/or near real-time analysis of even large-scale power systems. Using the modified IEEE 39 bus network, results show that the method can effectively identify stable OCs without the need to allocate more inertia resources. Nevertheless, an interesting area for future research work would be the application of the proposed method on realworld power networks which are larger and more complex.